\documentstyle{article}

\def\proof{\noindent{\bf Proof: }}
\def\qed{\hfill\mbox{$\Box$}\medskip}

\def\B{{\bf B}}

\def\R{{\bf R}}
\def\X{{\bf X}}
\def\Y{{\bf Y}}

\newcounter{hour}     
\newcounter{hours}    
\newcounter{minute}   
\newcounter{sixty}    

\newcommand{\setclock} {
\setcounter{hour}{\time}
\setcounter{sixty}{60}
\setcounter{minute}{\time}
\divide \value{hour} by \value{sixty}
\setcounter{hours}{\value{hour}}
\multiply \value{hours} by \value{sixty}
\advance \value{minute} by -\value{hours}
}

\newcommand{\clock}{
\the\value{hour}:\ifnum\the\value{minute}<10 0\fi\the\value{minute}
}

\setclock

\newtheorem{theorem}{Theorem}
\newtheorem{lemma}[theorem]{Lemma}
\newtheorem{cor}[theorem]{Corollary}

\title{Richman Games}

\author{Andrew J.~Lazarus\thanks{email: {\tt drlaz@aol.com} .} \\
Berkeley, California 
\and Daniel E. Loeb\thanks{email: {\tt loeb@labri.u-bordeaux.fr} . 
URL: {\tt http://www.labri.u-bordeaux.fr/\~{ }loeb/ }.
Partially supported by URA CNRS 1304,
EC grant CHRX-CT93-0400, the PRC Maths-Info, and NATO CRG 930554.}\\
Universit\'e de Bordeaux I
\and James G.~Propp\thanks{email: {\tt propp@math.mit.edu} .}
\\ The Massachusetts Institute of Technology
\and Daniel Ullman\thanks{email: {\tt dullman@math.gwu.edu} .} \\ The
George Washington University} 

\date{\today}

\begin{document}

\maketitle
\begin{center} Dedicated to David Richman, 1956--1991 
\end{center}

\vfill

\begin{abstract}
A Richman game is a combinatorial game in which, rather than alternating
moves, the two players bid for the privilege of making the next move.
We find optimal strategies for both 
the case where a player knows how much money his or her opponent has 
and the case where the player does not.					
\end{abstract}

\vfill

\bigbreak 


\newpage

\section{Introduction}\label{Intro}
There are two game theories.  The first is now sometimes referred to as
matrix game theory and is the subject of the famous von Neumann and
Morgenstern treatise \cite{vNM}.  In matrix games, two 
players make simultaneous moves and a payment is made from one player
to the other depending on the chosen moves.  
Optimal strategies often involve randomness and concealment of information. 

The other game theory is the combinatorial theory of {\em Winning Ways}
\cite{WW}, with origins back in the work of Sprague \cite{S} and 
Grundy \cite{G}, and largely expanded upon by Conway
\cite{ONAG}.  In combinatorial games, two players move alternately.
We may assume that each move consists of  
sliding a token from one vertex to another
along an arc in a directed graph.					
A player who cannot move loses.  There is no hidden information and 
there exist deterministic optimal strategies.

In the late 1980s,
David Richman suggested a class of games which share some aspects
of both sorts of game theory.  Here is the set-up:  The game is played
by two players (Mr.\ Blue and Ms.\ Red), each of whom has some money.
There is an underlying combinatorial game in which
a token rests on a vertex of some finite directed graph.
There are two special vertices, denoted by $b$ and $r$;
Blue's goal is to bring the token to $b$
and Red's goal is to bring the token to $r$.
The two players repeatedly bid for the right to make the next move.  
One way to execute this bidding process
is for each player to write secretly on a card a nonnegative real number 
no larger than the number of dollars he or she has; 
the two cards are then revealed simultaneously.  
Whoever bids higher pays the amount of the bid 
to the opponent and
moves the token from the vertex it currently occupies
along an arc of the directed graph to a successor vertex.
Should the two bids be equal, the tie is broken by a toss of a coin.  
The game ends when one player moves the token 
to one of the distinguished vertices.  
The sole objective of each player is to have the game end with the
token on his or her vertex. 
(At the game's end, money loses all value.)  
The game is a draw if neither distinguished vertex is ever reached.	

Note that with these rules (compare with \cite{Berl}), there is never a
reason for a negative bid: since all successor vertices are available to 
both players, it cannot be preferable to have the opponent move next. 
That is to say, there is no reason to part with money for the chance
that your opponent will carry out through negligence
a move that you yourself could perform through astuteness.

A winning strategy is a policy for bidding and moving that guarantees a
player the victory, given fixed initial data. (These initial data include 
where the token is, how much money the player has, and possibly how much 
money the player's opponent has.) 					
In section \ref{How}, we
describe a winning strategy for Richman games. In particular,
we prove the following facts, which might seem surprising:
\begin{itemize}
\item There exists a critical ratio $R(v)$ such that Blue
(respectively, Red) has a winning strategy from a vertex $v$
if Blue's share of the money, 
expressed as a fraction of the total money supply, 
is greater than (resp., less than) $R(v)$.
(This is not so surprising in the case of acyclic games,
but for games in general, one might have supposed it possible that,
for a whole range of initial conditions,
play might go on forever.)
\item There exists a strategy such that
if a player has more than~$R(v)$ and applies the strategy,
the player will win with probability 1,
without needing to know how much money the opponent has.		
\end{itemize}

In proving these assertions, it will emerge that an optimal
bid for Blue is $R(v)-R(u)$ times the total money supply,
where $v$ is the current vertex and $u$ is the successor
of $v$ for which $R(u)$ is as small as possible.  A player
who cannot bid this amount ``has already lost,'' in the sense
that there is no winning strategy for that player.  On the
other hand, a player who has a winning strategy of any kind
and bids $R(v)-R(u)$ will still have a winning strategy one
move later, regardless of who wins the bid, as long as he
or she is careful to move to $u$ if he or she does win the bid.

It follows that we may think of $R(v)-R(u)$ as the ``fair price'' 
that Blue should be willing to pay for the privilege of trading
the position $v$ for the position $u$.  Thus we may define
$1-R(v)$ as the {\it Richman value} of the position $v$, so
that the fair price of a move exactly equals the difference in
values of the two positions.  However, it is more convenient
to work with $R(v)$ than with $1-R(v)$, so we have opted for
this approach here.  We call $R(v)$ the {\it Richman cost} of
the position $v$.

We will see that for all $v$ other than the distinguished vertices
$b$ and $r$, $R(v)$ is the average of $R(u)$ and $R(w)$, where
$u$ and $w$ are the successors of $v$ in the digraph that minimize
and maximize $R(\cdot)$, respectively.  In the case where the digraph 
underlying the game is acyclic, this averaging-property makes it
easy to compute the Richman costs of all the positions, beginning 
with the positions $b$ and $r$ and working backwards.  On the other
hand, if the digraph contains cycles it is not so easy to work
out precise Richman costs.						

We defer most of our examples to another paper \cite{next} in which
we also consider infinite digraphs and discuss the
complexity of the computation of Richman costs.


\section{The Richman Cost Function}\label{How}				

Henceforth, $D$ will denote a finite directed graph $(V,E)$ with a	
distinguished blue vertex $b$ and a distinguished red vertex $r$ such
that from every vertex there is a path to at least one of the distinguished
vertices. For $v\in V$, let $S(v)$ denote the set of successors of $v$ in
$D$, that is $S(v)=\{w\in V \colon (v,w)\in E\}$. Given any function
$f\colon V\to 
[0,1]$, we define 
$$f^+(v)= \displaystyle \max_{w\in S(v)} f(w) \quad \mbox{and} \quad
f^-(v)= \min_{u\in S(v)} f(u).$$

The key to playing the Richman game on $D$ is to 
attribute costs to the vertices of $D$
such that the cost of every vertex (except the two distinguished
vertices) is 
the average of the lowest and highest costs of its successors. Thus, a
function $R\colon V\to [0,1]$ is called a {\em Richman cost function} 
if $R(b)=0$, $R(r)=1$, and for every other $v\in V$ we have 
$R(v) = (R^+(v) + R^-(v))/2$. 
(Note that Richman costs are a curious sort of variant
on harmonic functions on Markov chains \cite{Woess} where
instead of averaging over {\it 
all} the successor-values, we average only over the two extreme values.)
The relations $R^+(v) \geq R(v) \geq R^-(v)$ and $R^+(v) + R^-(v) = 2R(v)$
will be much used in what follows.
\begin{theorem}
\label{existence} 
The digraph $D$ has a Richman cost function~$R(v)$.
\end{theorem}
\proof 
We introduce an auxiliary function $R(v,t)$
whose game-theoretic significance 
will be made clearer in Theorem \ref{govern}.				
Let $R(b,t)=0$ and $R(r,t)= 1$ for all $t\in{\bf N}$. For $v\notin \{b,r\}$,
define $R(v,0)=1$ and $R(v,t)= (R^+(v,t-1) + R^-(v,t-1))/2$ for $t>0$. 
It is easy to see that $R(v,1) \leq R(v,0)$ for all $v$,
and a simple induction shows that $R(v,t+1) \leq R(v,t)$ for all $v$
and all $t \geq 0$.							
Therefore
$R(v,t)$ is weakly decreasing and bounded below by
zero as $t\to\infty$, hence convergent. It is also evident that the function
$\displaystyle v\mapsto\lim_{t\to\infty}R(v,t)$ satisfies the
definition of a Richman cost function. \qed

\paragraph{Alternate proof:} 
Identify functions $f \colon V(D)\to [0,1]$ with 
points in the $|V(D) |$-dimensional cube $Q=[0,1]^{|V(D)|}$.  
Given $f\in Q$, define $g\in Q$
by $g(b)=0$, $g(r)=1$, and, for every other $v\in V$,
$g(v) = (f^+(v) + f^-(v))/2$.  The map $f\mapsto g$ is clearly a
continuous map from $Q$ into $Q$, and so by the Brouwer fixed point
theorem it has a fixed point.  This fixed point is a Richman cost function. 
\qed

This Richman cost function does indeed govern the winning strategy,
as we now prove.
\begin{theorem}\label{govern} \label{inductlemma}
Suppose Blue and Red play the Richman game on the digraph $D$ 
with the token initially located at vertex $v$.
\begin{enumerate}
\item[{\bf (1)}] If Blue's share of the total money supply exceeds 
$\displaystyle R(v)=\lim_{t\rightarrow\infty } R(v,t) $, 
then he has a winning strategy.
\item[{\bf (2)}]  Moreover, his victory will require
at most $t$ moves if his share of the money supply exceeds $R(v,t)$. 
\end{enumerate}
\end{theorem}
\proof Without loss of generality, money may be scaled so that the 
total supply is one dollar. 
Whenever Blue has
over $R(v)$ dollars, he must have over $R(v,t)$ dollars for some $t$. 
We prove {\bf (2)} by induction on $t$. At $t=0$, Blue has over
$R(v,0)$ dollars 
only if $v=b$, in which case he has already won. 

Now assume {\bf (2)} is true for $t-1$, and let Blue have more
than $R(v,t)$ dollars.  There  
exist neighbors $u$ and $w$ of $v$ such that $R(u,t-1)=R^-(v,t-1)$ and 
$R(w,t-1)=R^+(v,t-1)$, so that $R(v,t) = (R(w,t-1)+R(u,t-1))/2$.
Blue can bid $(R(w,t-1)-R(u,t-1))/2$ dollars. If Blue wins the bid at $v$,
then he moves to~$w$ and 
forces a win in at most $t-1$ moves (by the induction hypothesis), since he 
has more than $(R(w,t-1)+R(u,t-1))/2-(R(w,t-1)-R(u,t-1))/2 = R(u,t-1)$ 
dollars left.  If Blue loses the bid, then Red will move to some $z$, but 
Blue now has over
$(R(w,t-1)+R(u,t-1))/2+(R(w,t-1)-R(u,t-1))/2 = R(w,t-1) \ge R(z,t-1)$ dollars,
and again wins by the induction hypothesis.
\qed									

One can define another auxiliary function $R'(v,t)$ where
$R'(b,t)=0$ and $R'(r,t)= 1$ for all $t\in{\bf N}$, and 
$R'(v,0)=0$ and $R'(v,t)= (R^+(v,t-1) + R^-(v,t-1))/2$ for $v\notin
\{b,r\}$, $t>0$.  By an argument similar to the proof of Theorem
\ref{govern}, this also converges to a Richman cost function $R'(v)\leq
R(v)$ (with $R(v)$ defined as in the proof of Theorem \ref{existence}).
Thus, $R'(v,t)$ indicates how much money Blue needs to prevent Red from
forcing a win from $v$ in $t$ or fewer moves, so $R'(v)$ indicates
how much money Blue needs to prevent Red from forcing a win in any length
of time.

For certain {\em infinite} digraphs, it can be shown \cite{next}
that $R'(v)$ is strictly less than $R(v)$. When Blue's share of the
money supply lies strictly between $R'(v)$ and $R(v)$, both players
can prevent the other player from winning. Thus, optimal play leads to a
draw. 

Nevertheless, in this paper, we assume that $D$ is finite, and we can
conclude that there is a unique Richman cost function $R'(v)=R(v)$. 
\begin{theorem}\label{uniqueness}
The Richman cost function of the digraph $D$ is unique.
\end{theorem}
The proof of Theorem~\ref{uniqueness} requires the following
definition and technical lemma.
An edge $(v,u)$ is said to be an {\em edge of steepest descent} if 
$R(u)=R^{-}(v)$.
\newcommand{\SD}[1]{\overline{#1}}
Let $\SD{v}$ be the transitive closure of $v$ under the 
steepest-descent relation.  That is, $w\in \SD{v}$ if there exists a path
$v=v_{0}, v_{1}, v_{2}, ..., v_{k}=w$
such that $(v_{i}, v_{i+1})$ is an edge of steepest descent
for $i=0,1,...,k-1$.							
\begin{lemma} 
\label{steeplemma}
Let $R$ be any Richman cost function of the digraph $D$. 
If $R(z) < 1$, then $\SD{z}$ contains~$b$.
\end{lemma}
\paragraph{Proof of Lemma:}  
Suppose $R(z)<1$.
Choose~$v\in \SD{z}$ such that
$\displaystyle R(v)=\min_{u\in 
\SD{z}} R(u)$. Such a~$v$ must exist because $D$ (and hence~$\SD{z}$)
is finite.  
If $v=b$, we're done.  Otherwise, assume $v \neq b$, and let $u$
be any successor of~$v$.  The definition of~$v$ implies $R^-(v)=R(v)$, 
which forces $R^+(v)=R(v)$. 
Since $R(u)$ lies between $R^-(v)$ and $R^+(v)$,
$R(u)=R(v)=R^-(v)$. Hence $(v,u)$ is an edge of steepest descent,
so $u \in \SD{z}$.  Moreover, $u$ satisfies the same defining property
that $v$ did (it minimized $R(\cdot)$ in the set $\SD{z}$), so the
same proof shows that for any successor~$w$ of~$u$, $R(w)=R(u)$
and $w\in\SD{z}$.  Repeating this, we see that
for any point~$w$ that may be reached from~$v$, $R(w)=R(v)$
and~$w\in\SD{z}$.  On the other hand, $R(r)$ is {\it not} equal to $R(v)$
(since $R(v) \leq R(z) < 1 = R(r)$),
so $r$ cannot be reached from $v$.
Therefore $b$ {\it can} be reached from $v$,
so we must have $b \in \SD{z}$. \qed					

\noindent {\bf Proof of Theorem~\ref{uniqueness}:}
Suppose that $R_1$ and $R_2$ are Richman cost functions of~$D$. Choose~$v$
such that $R_1-R_2$ is maximized at~$v$; such a~$v$ exists since~$D$
is finite. Let~$M=R_1(v)-R_2(v)$.
Choose $u_1, w_1, u_2, w_2$ (all successors of $v$)
such that $R^-_i(v)=R(u_i)$ and 
$R^+_i(v)=R(w_i)$. Since $R_1(u_1)\le R_1(u_2)$,
\begin{equation} \label{ineq:1}
R_1(u_1)-R_2(u_2) \le R_1(u_2)-R_2(u_2) \le M.
\end{equation}
(The latter inequality follows from the definition of~$M$.)
Similarly, $R_2(w_2) \ge R_2(w_1)$, so
\begin{equation}
\label{ineq:2}
R_1(w_1)-R_2(w_2) \le R_1(w_1)-R_2(w_1) \le M
\end{equation}
Adding~(\ref{ineq:1}) and~(\ref{ineq:2}), we have
$$(R_1(u_1)+R_1(w_1)) - (R_2(u_2)+R_2(w_2)) \le 2M.$$ The left side is
$2R_1(v)-2R_2(v)=2M$ so equality must hold in~(\ref{ineq:1}). In particular,
$R_1(u_2)-R_2(u_2) = M$; i.e., $u_2$ satisfies the hypothesis on~$v$.
Since $u_2$ was any vertex with
$R_2(u_2)=R^-_2(v)$, induction shows that $R_1(u)-R_2(u)=M$ for all $u\in
\SD{v}$ where descent is measured with respect to $R_2$. Since
$R_1(b)-R_2(b)=0-0=0$ and 
$b\in\SD{v}$, $R_1(v)-R_2(v) \leq 0$ everywhere. 
That is, $R_1 \leq R_2$.  The same argument for
$R_2-R_1$ shows the opposite inequality, so $R_1=R_2$.
\qed									

The uniqueness of the Richman cost function
implies in particular that the function $R'$
defined after the proof of Theorem~\ref{govern}~coincides
with the function $R$ constructed in the first proof of Theorem~\ref{existence}.
{}From this we deduce the following:
\begin{cor}\label{converse}
Suppose Blue and Red play the Richman game on the digraph $D$ 
with the token initially located at vertex $v$.
If Blue's share of the total money supply is less than
$\displaystyle R(v)=\lim_{t\rightarrow\infty } R(v,t) $, 
then Red has a winning strategy.\qed 
\end{cor}
It is also possible to reverse the order of proof, and to derive
Theorem \ref{uniqueness} from Corollary \ref{converse}.  For,
if there were two Richman functions $R_1$ and $R_2$,
with $R_1(v) < R_2(v)$, say, then by taking a situation in
which Blue's share of the money was strictly between $R_1(v)$
and $R_2(v)$, we would find that both Blue and Red had winning
strategies, which is clearly absurd.					

Theorem~\ref{govern} and Corollary \ref{converse}
do not cover the critical case where Blue has exactly $R(v)$ 
dollars. In the critical case, with both players using optimal strategy,
the outcome of the game depends on the outcomes of the coin-tosses used
to resolve tied bids.  Note, however, that in all other cases, the 
deterministic strategy outlined in the proof of Theorem~\ref{inductlemma}
works even if the player with the winning strategy concedes all ties
and reveals his intended bid and intended move before the bidding. 	

Summarizing Theorem \ref{govern} and Corollary \ref{converse},
we may say that								
if Blue's share of the total 
money supply is less (respectively, greater) than $R(v)$, then Red
(resp., Blue) has a winning strategy. 


\section{Other Interpretations}

Suppose the right to move the token is 
decided on each turn by the toss of a fair coin.
Then induction on $t$ shows that the probability
that Red can win from the position $v$ in at most $t$ moves
is equal to $R(v,t)$, as defined in the previous section.
Taking $t$ to infinity, we see that $R(v)$ is equal to
the probability that Red can force a win against optimal play by Blue.
That is to say, if both players play optimally,
$R(v)$ is the chance that Red will win.
The uniqueness of the Richman cost function tells us that
$1-R(v)$ must be the chance that Blue will win.
The probability of a draw is therefore zero.

If we further stipulate that the moves themselves must be random,
in the sense that the player whose turn it is to move
must choose uniformly at random from among the finitely many legal options,
then we do not really have a game-like situation anymore;
rather, we are performing a random walk on a directed graph
with two absorbing vertices,
and we are trying to determine the respective probabilities of absorption
at these two vertices.
In this case, the relevant probability function
is just the harmonic function on the digraph $D$ 
(or, more properly speaking, the harmonic function
for the associated Markov chain \cite{Woess}).

Another interpretation of the Richman cost,
brought to our attention by Noam Elkies, comes from a problem
about makings bets.  Suppose you wish to bet (at even odds)
that a certain baseball team will win the World Series, but that
your bookie only lets you make even-odds bets on the 
outcomes of individual games.  Here we assume that the winner of 
a World Series is the first of two teams to win four games.  To
analyze this game, we create a directed graph whose vertices
correspond to the different possible combinations of cumulative
scores in a World Series, with two special terminal vertices
(blue and red) corresponding to victory for the two respective teams.  
Assume that your initial amount of money is $1/2$ a grand (\$500),
and that you want to end up with either $0$ grand or $1$ grand,
according to whether the blue team or the red team wins the Series.
Then it is easy to see that the Richman cost at a vertex tells 
exactly how much money you want to have left if the
corresponding state of affairs transpires, and that the amount
you should bet on any particular game is the common value of
$R(v)-R(u)$ and $R(w)-R(v)$, where $v$ is the current position,
$u$ is the successor position in which Blue wins the next game,
and $v$ is the successor position in which Red wins the next game.

\section{Incomplete Knowledge}

Surprisingly, knowledge of one's opponent's money supply is
unnecessary for the construction of a winning strategy.

Define Blue's {\em safety ratio} at $v$ to
be the fraction of the total money that he has in his possession,
divided by $R(v)$ (the fraction that he needs in order to win).
Note that Blue will not know the value of his safety ratio,
since we are assuming that he has no idea
how much money Red has.

\begin{theorem} \label{thm2}
Suppose Blue has a safety ratio strictly greater than 1.
Then Blue has a strategy that wins with probability~$1$ and
does not require knowledge of Red's money supply.  
If, moreover, the digraph $D$ 
is acyclic, then his strategy wins regardless of tiebreaks;
that is, ``with probability 1'' can be replaced by ``definitely''.	
\end{theorem}
\newcommand{\Rcrit}{\R_{\rm crit}}
\newcommand{\Ractual}{\R_{\rm actual}}
%
\proof 
Here is Blue's strategy: When the token is at vertex~$v$,
and he has $\B$~dollars, he should act as if his safety ratio is
1; i.e., he should play as if Red has
$\Rcrit$ dollars with $\B/(\B+\Rcrit) = R(v)$ and 
the total amount of money is $\B + \Rcrit = \B/R(v)$ dollars. 
He should accordingly bid $$\X=\frac{R(v)-R^-(v)}{R(v)}\B$$ dollars.
Suppose Blue wins (by outbidding or by tiebreak) and moves to $u$ along an
edge of steepest descent.  Then Blue's safety ratio changes 
$$\mbox{from \ \ \ } \frac{\left( \frac{ \B }{\B+\R} \right)}{R(v)}
\ \ \ \mbox{to   \ \ \ } \frac{\left( \frac{\B-\X}{\B+\R} \right)}{R(u)},$$
where $\R$ is the actual amount of money that Red has.
However, these two safety ratios are actually equal, since 
$$\frac{\B-\X}{\B} = 1 - \frac{\X}{\B} = 
1 - \frac{R(v)-R(u)}{R(v)} = \frac{R(u)}{R(v)}.$$
Now suppose instead that Red wins the bid (by outbidding or by tiebreak)
and moves to $z$.  Then Blue's safety ratio changes
$$\mbox{from \ \ \ } \frac{\left( \frac{ \B }{\B+\R} \right)}{R(v)}
\ \ \ \mbox{to   \ \ \ } \frac{\left( \frac{\B+\Y}{\B+\R} \right)}{R(z)},$$
with $\Y \geq \X$.  Note that the new safety ratio is greater than or equal to
$$\frac{\left( \frac{\B+\X}{\B+\R} \right)}{R(w)},$$
where $R(w)=R^+(v)$.  But this lower bound on the new safety ratio
is equal to the old safety ratio, since
$$\frac{\B+\X}{\B} = 1 + \frac{\X}{\B} = 
1 + \frac{R(w)-R(v)}{R(v)} = \frac{R(w)}{R(v)}.$$

In either case, the safety ratio is non-decreasing.
In particular, the safety ratio must stay greater than 1.
On the other hand, if Blue ever loses,
then his safety ratio at that moment would have to be at most 1,
since his fraction of the total money supply
cannot be greater than $R(r)=1$.
Consequently, our assumption that Blue's safety ratio started out
being greater than 1 
implies that Blue can never lose.
In an acyclic digraph, infinite play is impossible, so
the game must terminate at~$b$ with a victory for Blue.

In the case where cycles are possible, 
suppose first that at some stage Red outbids Blue by $\epsilon\B > 0$
and gets to make the next move,
say from $v$ to $w$.
If Blue was in a favorable situation at $v$,
then the total amount of money that the two players have between them
must be less than $\frac{1}{R(v)}\B$.
On the other hand, 
after the payoff by Red, Blue has
$\B+\X+\epsilon\B=(1+\frac{R(v)-R(u)}{R(v)}+\epsilon)\B
=(\frac{2R(v)-R(u)}{R(v)}+\epsilon)\B
=(\frac{R(w)}{R(v)}+\epsilon)\B$,
so that Blue's total share of the money must be more than
$R(w)+\epsilon R(v)$.
Blue can do this calculation as well as we can;
he then knows that if he had been in a winning position to begin with, his
current share of the total money must exceed $R(w)+\epsilon R(v)$.
Now, $R(w)+\epsilon R(v)$ is
greater than $R(w,t)$ for some $t$, so Blue can win in $t$ moves.
Thus, Red loses to Blue's
strategy if she ever bids more than he does.
Hence, if she hopes to avoid losing,
she must rely entirely on tiebreaking. 
Since Blue cannot lose, there is a shortest path of steepest descent 
to~$b$ from every vertex that will be reached as play continues. 
Let the longest of these paths have length~$N$. 
Then Blue will win the game when he wins~$N$ consecutive tiebreaks 
(if not earlier).
\qed									

\begin{figure}[htb]
\caption{The digraph $D$ and its Richman costs}
\vskip 0.3in 
\label{cycles}
\setlength{\unitlength}{.2cm} 
\begin{center}
\begin{picture}(40,15)(-15,-10)
\put(10,0){\circle{1.5}}
\put(0,0){\circle{1.5}}
\put(20,5){\circle{1.5}}
\put(20,-5){\circle{1.5}}
\put(-10,5){\circle{1.5}}
\put(-10,-5){\circle{1.5}}
\put(1,0){\vector(1,0){8}}
\put(11,.5){\vector(2,1){8}}
\put(11,-.5){\vector(2,-1){8}}
\put(-9,4.5){\vector(2,-1){8}}
\put(-1.5,-0.5){\vector(-2,-1){8}}
\put(-10,-4){\vector(0,1){8}}
\put(20,5){\makebox(0,0){$b$}}
\put(20,-5){\makebox(0,0){$r$}}
\put(0,0){\makebox(0,0){$v$}}
\put(10,2){\makebox(0,0){$1/2$}}
\put(0,2){\makebox(0,0){$1/2$}}
\put(20,7){\makebox(0,0){$0$}}
\put(20,-3){\makebox(0,0){$1$}}
\put(-10,7){\makebox(0,0){$1/2$}}
\put(-12,-5){\makebox(0,0){$1/2$}}
\end{picture}
\end{center}
\end{figure}
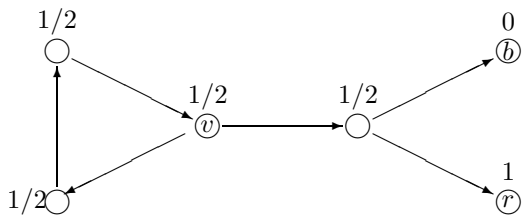

When $D$ has cycles, Blue may need to rely on tiebreaks
in order to win, as in the case of the Richman game played on the
digraph pictured in Figure~\ref{cycles}.
Suppose that the token is at vertex marked $v$, that
Blue has $\B$ dollars, and that Red has $\R$ dollars.  
Clearly, Blue knows he can win the game if $\B > \R$.
But without knowing $\R$, it would be
imprudent for him to bid any positive amount $\epsilon\B$
for fear that Red actually started with $(1-\epsilon)\B$ dollars;
for if that were the case, and his bid were to prevail,
the token
would move to a vertex where the Richman cost
is $1/2$ and Blue would have less money than
Red. Such a situation will lead to a win for Red if
she follows the strategy outlined in Theorem~\ref{inductlemma}. 

\section{Rationality}

For every vertex $v$ of the digraph $D$ (other than $b,r$),
let $v^+$ and $v^-$ denote successors of $v$
for which $R(v^+)=R^+(v)$ and $R(v^-)=R^-(v)$.
Then we have $R(b) = 0$, $R(r) = 1$, and
$2R(v) = R(v^+) + R(v^-)$ for $v \neq b,r$.
We can view this as a linear program.
By Theorem \ref{uniqueness}, this system must have a unique solution.
Since all coefficients are rational,
we see that Richman costs are always rational numbers.

The linear programming approach also gives us a conceptually simple
(though computationally dreadful) way to calculate Richman costs.
If we augment our program by adding additional conditions of the form
$R(v^-) \leq R(w)$ and $R(v^+) \geq R(w)$
where $v$ ranges over the vertices of $D$ other than $b$ and $r$
and where, for each $v$,
$w$ ranges over all the successors of $v$,
then we are effectively adding in the constraint
that the edges from $v$ to $v^-$ and $v^+$
are indeed edges of steepest descent and ascent, respectively.
The uniqueness of Richman costs tells us that
if we let the mappings $v \mapsto v^-$ and $v \mapsto v^+$
range over all possibilities
(subject to the constraint that both $v^-$ and $v^+$
must be successors of $v$),
the resulting linear programs (which typically will have no solutions at all)
will have only solutions that correspond to the genuine Richman cost functions.
Hence, in theory one could try all the finitely many possibilities
for $v \mapsto v^-$ and $v \mapsto v^+$ 
and solve the associated linear programs
until one found one with a solution.
However, the amount of time such an approach would take
increases exponentially with the size of the directed graph.
In \cite{next}, we discuss more efficient approaches.			

We will also discuss in \cite{next}, among other things,
a variant of Richman games,
which we call ``Poorman'' games,
in which the winning bid is paid to a third party
(auctioneer or bank) rather than to one's opponent.
The whole theory carries through largely unchanged,
except that Poorman costs are typically irrational.

\end{document}